\documentclass[preprint,12pt]{elsarticle}

\usepackage{amsmath, amsxtra, amsfonts,amscd,amssymb}
\usepackage{graphicx,wrapfig}
\usepackage{epstopdf}
\usepackage{url}
\usepackage[algo2e,linesnumbered, vlined,ruled]{algorithm2e}
\usepackage{float}
\usepackage{multirow}

\def\A{\mathcal{A}}

\newtheorem{theorem}{Theorem}[section]

\newtheorem{definition}[theorem]{Definition}

\begin{document}

\begin{frontmatter}



\title{Optimal Impulse Control  of a Simple Reparable  System in a Nonreflexive Banach Space}


\author{Weiwei Hu$^a$, Rongjie Lai$^b$, Houbao Xu$^c$, Chuang Zheng$^d$}
%

\address{$^a$Department of Mathematics, Oklahoma State University, Stillwater, USA. weiwei.hu@okstate.edu }
\address{$^b$Department of Mathematics, Rensselaer Polytechnic Institute, Troy, USA.  lair@rpi.edu }
\address{$^c$Department of Mathematics, Beijing Institute of  Technology, Beijing, China.  xuhoubao@bit.edu.cn}
\address{$^d$School of Mathematical Science, Beijing Normal University, Beijing, China.  chuang.zheng@bnu.edu.cn}
\begin{abstract}
We  discuss the problem of  optimal impulse  control representing the preventive maintenance of a simple  reparable system.  The  system model  is  governed by  coupled transport and integro-differential equations in a nonreflexive Banach space.  The objective of this paper is to construct  nonnegative impulse control inputs at  given system running times  that  minimize the probability of the system in  failure mode.  To guarantee the  nonnegativity  of the controlled system, we consider the control inputs to depend on the system state.  This essentially leads to a  bilinear control problem.  We first present a rigorous proof of existence of an optimal controller and  then apply the variational inequality  to derive the first-order necessary conditions of optimality.

\end{abstract}


\end{frontmatter}



\section{Introduction}
\label{sec:intro}
Reparable systems  occur naturally in many real  world problems   such as  product design, inventory systems, computer networking and complex manufacturing processes, etc. In recent years,  mathematical models governed by distributed parameter systems of coupled   partial and ordinary hybrid equations have been widely used to study the  reliability of  reparable systems (cf.~\cite{Chung-1, DHR-1, Gupur-1, HG-1, HR-1, HR-2, HXYZ-1, XYZ-1}).  Reliability  is defined as the probability that the system, subsystem or component will operate successfully by a given time $t$ (cf.~\cite{bazovsky2004reliability, Sandler2012}). This paper is mainly concerned  with the optimal control problem of  a reparable multi-state system  introduced by Chung \cite{Chung-1},
 which is described as coupled transport and integro-differential equations. In particular, we consider that the system has one  failure state in our current work. This will be sufficient to capture the essence  of the control design. The precise model of system equations is described  by
\begin{align}
&\frac{dp_0(t)}{d t} =-\lambda_1  p_0(t)+ \int_0^{l} \mu_1(x)p_1(x,t) \,d x,  \label{eq1}   \\
&\frac{\partial{p_1(x,t)}}{\partial{t}} + \frac{\partial{p_1(x,t)}}{\partial{x}} =-\mu_1(x)p_{1}(x,t),
 \label{eq2}
\end{align}
with the boundary condition
\begin{equation}
p_1(0,t)=\lambda_1 p_0(t),
\label{BC}
\end{equation}
and initial conditions
\begin{equation}
p_0(0)=1, \quad p_1(x,0)=0.
\label{IC}
\end{equation}
Here 
$p_{0}(t)$ stands for the probability that the device is in good state, represented as $0$, at time $t$; $p_{1} (x, t)$ stands for
the probability density (with respect to the repair time) that the failed device is in  failure state, represented by $1$, and has an elapsed repair time of $x$ at time $t$, where $x\in [0, l]$ with $l<\infty$; $\lambda_{1}$ stands for the
constant failure rate of the system for failure mode $1$; $\mu_{1}(x)$ stands for the
time-dependent repair rate when the device is in state $1$ and has an elapsed repair time of $x$;
probability $\hat{p}_{1}(t)$ that the failed device is in state $1$ at time $t$ is defined by
   \begin{align}
            \hat{p}_1(t)=\int^{l}_0p_1(x,t)\,d x, \quad t\geq 0. \label{prob_p1}
    \end{align}
  The following assumptions are associated with the device:
\begin{enumerate}
\item The failure rates are constant;
\item  All failures are statistically independent;
\item  All repair times of failed devices are arbitrarily distributed;
\item The repair process begins soon after the device is in failure state;
\item The repaired device is as good as new;
\item  No further failure can occur when the device has been down.
\end{enumerate}

Furthermore, we assume that the repair rate has the following properties
  \begin{equation}
\int^{l'}_{0} \mu_{1}(x)\,d x <\infty, \quad l'<l, 
\quad \text{and}\quad
 \int^{l}_{0} \mu_{1}(x)\,d x =\infty.
 \label{Assp_repair_fun1}
\end{equation}

The pointwise and steady-state availability of the uncontrolled multi-state system was discussed in \cite{Chung-1} by solving the inversion of the Laplace transformation. However, this approach used two potential assumptions that the system has a nonnegative time-dependent solution and the solution is asymptotically stable, which are nontrivial when the 
repair rate is time dependent.  Xu, Yu and Zhu in \cite{XYZ-1} provided  a rigorous mathematical framework for proving the well-posedness and asymptotic stability of the  system by using $C_0$-semigroup theory. It is proved that the system operator generates a positive
$C_{0}$-semigroup of contraction, therefore the system has a unique nonnegative time-dependent solution.
It is also shown  that 0 is a simple eigenvalue of the system operator and the unique
spectral point on the imaginary axis. In particular, the system is conservative in the sense that the sum of  the probabilities of the system in good mode and failure mode is always $1$. Moreover, Hu, Xu, Yu and Zhu in \cite{HXYZ-1} 
showed that  the $C_{0}$-semigroup is quasi-compact and irreducible. As a result,  it follows that the
time-dependent solution converges to the steady-state solution  exponentially, which is the positive
eigenfuction corresponding to the simple eigenvalue $0$. However, the previous work is based on the assumption that the repair time $l=\infty$,
which is not realistic.  Neither the system  can be under repair nor the server can work forever.  Recently,  Hu in  \cite{hu2016differentiability}  considered $l<\infty$  and further proved that the $C_0$-semigroup is eventually compact and  eventually differentiable.  In this case, condition \eqref{Assp_repair_fun1} indicates that if the system, subsystem or component  can not be repaired over the finite repair time interval  $[0, l]$, then it will be replaced by the new one immediately. Furthermore, the problems of  controllability and enhancement of stabilizability of this system   by distributed controls  have been discussed in \cite{huAcc2017, WeiHuZheng}.

This paper aims at minimizing  the probability of the system at failure mode by employing a maintenance policy, which is interpreted as the control inputs.    Maintenance is defined as any action that restores failed units to an operational condition or retains non-failed units in an operational state (cf.~\cite{jardine2013maintenance, moubray1997reliability}).  These actions affect  the overall performance of the system such as reliability, availability,
downtime, cost of operation, etc. 
A proper  maintenance policy  and a feasible approach are crucial.  The
optimal control design formulated  in this paper provides insight into the development of such a policy.  In general, there are three types of maintenance actions: corrective maintenance, preventive maintenance and inspections.  Corrective maintenance serves to restore a failed system
to operational status. As we can see from the model equations \eqref{eq1}--\eqref{eq2},  repair rate $\mu_1(x)$ plays the role of corrective
maintenance. Since a
component's failure time is not known a priori, $\mu_1(x)$ is performed during unpredictable intervals. This usually involves replacing or repairing the component that is
responsible for the failure of the overall system. 

Compared to corrective maintenance, the concept of preventive maintenance is  to replace components or subsystems before they fail in order to promote
 continuous system operation. Cost needs to be taken into account  in preventive maintenance since financially it is more sensible to replace parts or
 components that have not failed at predetermined intervals rather than to wait for a system failure. The latter  may result in a costly disruption in operations.   In our current work, we propose a  preventive maintenance policy that are applied   at the given system running times. To formulate a meaningful strategy, two criteria need to be satisfied in our control design:
 \begin{enumerate}
 \item   The controlled system is  nonnegative and conservative;
 \item  The control inputs are nonnegative.
 \end{enumerate}

The paper is organized as follows. In section 2, we introduce an impulse control design to the reparable system and establish the well-posedness of the controlled system in a nonreflexive Banach space.    Then we prove the existence of an optimal control in section 3.  In section 4,  we present the first-order necessary conditions of optimality by using a  variational inequality. 
 \section{Optimal Impulse Control Design}
 \label{sec:OptControl}
Consider that the reparable system (\ref{eq1})--(\ref{IC}) is  controlled by  the impulse control inputs at  given system running times $t_{i}, 0\leq t_{1}<t_{2}<\dots<t_{N}\leq T$,  for some  $T>0$. 
Let  $u_{i}(x)\geq 0,  x\in [0, l], i=1,2, \dots, N,$ be the corresponding input  intensities, which represent the  the update/replacement rate when the system  is in failure state at time $ t_{i}$ with an elapsed repair time $x$. 
The controlled system  is described by
 \begin{align}
&\frac{dp_0(t)}{d t}=-\lambda_1  p_0(t)+ \int_0^{l} \mu_1(x)p_1(x,t) \, d x 
+ \sum^{N}_{i=1}  \delta(t-t_{i}) \int_0^{l}u_{i}(x)\,d x, \label{eq_controlled1}\\
&\frac{\partial{p_1(x,t)}}{\partial{t}}+\frac{\partial{p_1(x,t)}}{\partial{x}} = -\mu_1(x)p_{1}(x,t) 
-\sum^{N}_{i=1} \delta(t-t_{i}) u_{i}(x),
\label{eq_controlled2}
\end{align}
where $\delta(t-t_{i})$  is the Dirac $\delta-$function supported at $t_{i}$.
Note  that the system  will not be updated/replaced if it is in good state. 
To study  the behavior of the reparable system in terms of probabilities, we consider  the state space $X=\mathbb{R} \times L^1[0,l] $ with $\|\cdot\|_X=|\cdot|+ \|\cdot\|_{L^1[0, l]}$, which is a nonreflexive Banach space.  Let $X^{*}$ be the dual of $X$, then $X^*=\mathbb{R}\times L^{\infty}[0, l]$. The duality between $X$ and $X^*$ is defined by
 \[(P, Q)=p_0q_0+\int^{l}_{0}p_1q_1\, dx, \]
for every $P=[p_0, p_1]^T\in X$ and  $Q=[q_0, q_1]^T\in X^*$.
The objective of this paper is to establish an optimal update/replacement policy  that minimizes the probability   $\hat{p}_1(t)$ of the system in  failure mode over  $[0, T]$, where  $\hat{p}_1(t)$ is defined by \eqref{prob_p1}.

To be physically meaningful,  we seek for  the control inputs $u_{i}(x)\geq 0$ with $x\in [0, l]$, such that  the controlled system is nonnegative and conservative. In particular,   we consider that $u_{i}$ depends on the probability density of the system in failure mode at $t_i$, that is,
\begin{eqnarray}
u_{i}(x)=b_{i}(x) p_1(x, t_i), \quad x\in [0, l], \quad i=1,2\dots, N,  \label{feedback}
\end{eqnarray}
where $0\leq b_{i}(x)\leq 1$, stands for the updated/replaced rate depending  on the elapsed repair time $x$.  In other  words, the input  intensity $u_{i}(x)$   at  time $t_{i}$ is up to the 
probability density  of the device in failure mode with an elapsed repair time of $x$ at time $t_{i}$.
For a given final time $T>0$, the current work   aims at deriving  the optimal distribution of $b_i(x)$ such that the  following cost functional is minimized:
\begin{align}
J( b_1, b_2, \dots, b_N)  = &\frac{1}{2} \int^{T}_{0}|\int^l_{0}p_1(x,t)\, dx|^2\,d t
+\frac{1}{2}\sum^N_{i=1}\beta_i \int^{l}_{0} |b_{i}(x)|^2 \,dx \nonumber\\
&+\frac{1}{2}| \int^l_{0}p_1(x,T)\, dx|^2,
\label{cost}
\end{align}
where $\beta_i\geq 0, i=1,2, \dots, N,$ are the  control weight  parameters. This essentially becomes  a bilinear optimal control problem. 

We first show that the controlled system is nonegative and conservative. Define the system operator $\A$ and its domain
 \begin{eqnarray}
\A P=
    \left[\begin{array}{c}
-\lambda_{1}  p_0+ \int_0^{l }\mu_1(x) p_1(x) \,d x\\
-(\frac{d }{d x}+\mu_1(x) )p_1(x)
\end{array}\right], \label{system_oper}
\end{eqnarray} 
\begin{eqnarray*}
D(\A)=\big\{P\in X\big| \ \frac{ d p_1(x)} {d x} \in L^1[0, l],  \int^{l}_{0}\mu_1(x)p_{1}(x)\, dx <\infty,  \ \text{and}\ p_1(0)=\lambda_1 p_0 \big\}.
 \end{eqnarray*}
Let $\mathcal{S}(t)$, $t\geq 0,$ denote the    $C_0$-semigroup generated by $A$, which is  a positive semigroup  of contraction and eventually compact and  eventually differentiable for $l<\infty$ \cite{hu2016differentiability, XYZ-1}.
Moreover, $0$ is a simple eigenvalue of $\A$ and  the only spectrum on the imaginary axis.    

Now let $f_{0}( t-t_{i})=\delta(t-t_{i})     \int^{l}_{0}b_{i}(x)p_{1}(x, t_{i})\,d x$, $f_{1}(x, t-t_{i})=- \delta(t-t_{i})  b_{i}(x)p_{1}(x, t_{i})$, and  $f(x,t-t_{i})=(f_{0}(t-t_{i}), f_{1}(x,t-t_{i}))^{T}$.
Then the controlled system \eqref{eq_controlled1}--\eqref{eq_controlled2} can be  rewritten as an abstract Cauchy  initial value problem in state space $X$
\begin{align}
\dot{P}(t)=\A P(t)+ \sum^{N}_{i=0} f(t-t_{i}), \quad t>0, \label{IVP_A}
\end{align}
with a general initial condition
\begin{align}
P(0)=(p_0(0), p_{1}(x, 0))^T \geq 0, \label{IVP_IC1}
\end{align}
satisfying
\begin{align}
p_0+\int^l_{0}p_{1}(x, 0)\, dx=1. \label{IVP_IC2}
\end{align}

The problem of impulse control has been discussed  in \cite{ahmed2011optimal} and the references cited therein. The solution of Cauchy problem \eqref{IVP_A}--\eqref{IVP_IC2} can be given by the  variation of parameter formula
\begin{align*}
P(t)=&\mathcal{S}(t)P(0)+ \sum^{N}_{i=1}\int^{t}_{t_{i}}\mathcal{S}(t-s)f(s-t_{i})\, d s\nonumber\\
=&\mathcal{S}(t)P(0)+\mathcal{S}(t-t_{1})\left[  \int^{l}_{0}b_{1}(x)p_{1}(x, t_{1})\, d x , - b_{1}(x)p_{1}(x, t_{1})\right]^{T}\nonumber\\
&+\mathcal{S}(t-t_{2})\left[  \int^{l}_{0}b_{2}(x)p_{1}(x, t_{2})\, d x , - b_{2}(x)p_{1}(x, t_{2})\right]^{T}\nonumber\\
&+\dots+\mathcal{S}(t-t_{N})\left[  \int^{l}_{0}b_{N}(x)p_{1}(x, t_{N})\, d x , - b_{N}(x)p_{N}(x, t_{N}))\right]^{T}. \nonumber\\
\label{Var_sol0}
\end{align*}
In particular, if  the control input is only exerted at $t=0$, then it becomes a start control problem. We  summarize the property of the solution in the following theorem.
\begin{theorem}
\label{T01}
For $t\in (t_{i}, t_{i+1}], 1\leq k\leq N-1$, the function given by the variation of parameter formula
\begin{equation}
P(t)=\mathcal{S}(t-t_{i})P(t_{i})+\mathcal{S}(t-t_{i}) \left[ \int^{l}_{0}b_{i}(x)p_{1}(x, t_{i})\, d x , - b_{i}(x)p_{1}(x, t_{i})\right]^{T},  \label{Var_sol1}
\end{equation}
is the  mild solution of the Cauchy initial value problem (\ref{IVP_A}).
Moreover, for $t\in (0,  t_{1}]$,
\begin{equation}
P(t)=\mathcal{S}(t)P(0), \label{Var_sol2}
\end{equation}
and for $t\in (t_{N}, \infty)$,
\begin{equation}
P(t)=\mathcal{S}(t-t_{N})P(t_{N}).  \label{Var_sol3}
\end{equation}
  \end{theorem} 
  Note that integrating  the second equation of \eqref{eq_controlled1} with respect to $x$ from $0$ to $l$, and then adding them to the first equation result in
\begin{eqnarray}
\frac{dp_{0}(t)}{dt}+\frac{d}{dt}\int^{l}_{0}p_{1}(x,t)\, d x=0. \label{conservation1}
\end{eqnarray}
This implies that 
\begin{eqnarray}
p_{0}(t)+\int^{l}_{0}p_{1}(x,t)\, d x=p_{0}(0)+ \int^{l}_{0}p_{1}(x, 0)\, d x=1, \quad \forall t\geq 0, \label{conservation2}
\end{eqnarray}
or
\begin{align}
\Vert P(t)\Vert_{X}=\Vert P(0)\Vert_{X}=1,  \quad \forall t\geq 0. \label{conservation3}
\end{align}
Thus the  system is conservative with respect to $\|\cdot\|_{X}$. In other words,  the sum of probability distributions of the controlled system  is always $1$ for every $t\geq 0$. 

Next replacing $p_{0}(t_{i})$ by $1-\int^{l}_{0}p_{1}(x, t_{i})\, d x$ in  (\ref{Var_sol1}), we have
\begin{align}
P(t)& =\mathcal{S}(t-t_{i})\left[p_{0}(t_{i})+ \int^{l}_{0}b_{i}(x)p_{1}(x, t_{i})\, d x , (1- b_{i}(x))p_{1}(x, t_{i})\right]^{T}\nonumber\\
 &=\mathcal{S}(t-t_{i})\left[1-\int^{l}_{0}(1-b_{i}(x))p_{1}(x, t_{i})\, d x , (1- b_{i}(x))p_{1}(x, t_{i})\right]^{T}, \nonumber\\ \label{P_theta_limit}
\end{align}
which is nonnegative for $0\leq b_{i}(x)\leq 1$.   
Further note that the solution has jumps at $t_{i}, i=1,\dots, N$. In fact, 
\begin{align}
P(t^{+}_{i})-P(t_{i})=&\mathcal{S}(t^{+}_{i}-t_{i})P(t_{i})\nonumber\\
&+\mathcal{S}(t^{+}_{i}-t_{i})\left[  \int^{l}_{0}b_{i}(x)p_{1}(x, t_{i})\, d x , - b_{i1}(x)p_{1}(x, t_{i})\right]^{T} 
-P(t_{i})\nonumber\\
=&\left[ \int^{l}_{0}b_{i}(x)p_{1}(x, t_{i})\, d x , - b_{i}(x)p_{1}(x, t_{i}))\right]^{T}.\label{sol_jump}
 \end{align}
Therefore, we have $P(t)\in PWC_l([0, T]; X)$, where $PWC_l([0, T]; X)$ denotes the space of piecewise continuous functions on $[0, T]$ with
values in $X$, that  are left continuous and possess righthand limits.

\section{Existence of an Optimal Solution}

 In this section, we address the existence of an optimal control  to problem \eqref{cost} subject to the controlled system \eqref{IVP_A}--\eqref{IVP_IC2}. Since the controlled system does not have any ``smoothing effects" on state $p_{1}(\cdot, t_i)$, we need an additional  condition imposed on $b_i$ to have the compactness of the admissible control set in order to handle the bilinear term $b_i(\cdot)p_{1}(\cdot, t_i)$.  
 To this end, we define the set of  admissible controls  by
 \begin{align}
U_{ad}= \{U=[b_1, b_2, \dots, b_N]^T\in (L^{2}[0, l])^N\big | &\  0\leq b_{i}\leq 1, i=1,2,\dots, N, \nonumber \\
& \text{is equicontinuous on}\ [0, l] \}.\label{U_ad}
\end{align}
Next we  introduce the  weak solution to \eqref{IVP_A}--\eqref{IVP_IC2}.
 \begin{definition} 
 \label{def1}
 For $P_0\in X$ with $\|P_0\|_{X}=1$, $P=[p_0, p_1]^T\in PWC_l([0, T]; X)$ is said to be a weak solution of system \eqref{IVP_A}--\eqref{IVP_IC2}, if $P$ satisfies 
 \begin{align}
( \dot{P}(t),Q)=(\A P(t), Q)+( \sum^{N}_{i=0} f(t-t_{i}), Q), \quad \forall Q=[q_0, q_1]^T\in X^*, \label{def_1}
 \end{align}
 and the initial condition 
 \begin{align}
 P(0)=P_0. \label{def_ini}
 \end{align}

 \end{definition}

 The following theorem provides the existence of an optimal solution. 
\begin{theorem}
\label{T04}
 There exists an optimal solution $(U^{*}, P^*)$ of problem \eqref{cost} subject to the controlled  system \eqref{IVP_A}--\eqref{IVP_IC2} in the sense of Definition \ref{def1}.
\end{theorem}
{\noindent {\textbf{Proof.} } 
Since $J$ is bounded from below, we can choose a minimizing sequence
  $\{U^{k}\}\subset U_{ad}$ for each $i=1, \dots, N$ such that
 \begin{align}
 \lim_{k\to \infty}J(U^{k})=\inf_{U\in U_{ad}}J(U).
 \end{align}
By the definition of $U_{ad}$, the sequence $\{U^{k}\}$ is uniformly bounded and  equicontinuous in $U_{ad}$. With the help of Ascoli's Theorem \cite{reed1980methods}, we may extract a subsequence, still denoted by $\{U^{k}\}$, such that 
\begin{align}
 U^{k}=[b^k_1, b^k_2, \dots, b^k_N]^T \to  U^*=[b^*_1, b^*_2, \dots, b^*_N]^T  \quad \text{uniformly  in}\quad [0, l]. \label{U_uniform}
 \end{align}
Let $P^k$ be the solutions of \eqref{def_1} corresponding to $U^k$ and satisfying the initial condition $P^k(0)=P_0$. 
Note  that $P^k\in PWC_{l}([0, T], X)$ and   $\|P^k(t)\|_{X}=1$ for any $t\geq 0$ based on \eqref{conservation2}. Thus for $t=t_i$, there exists a subsequence,  still denoted by $\{P^k\}$, satisfying 
\begin{align}
 P^{k}(\cdot, t_i)\to  P^*(\cdot, t_i) \quad \text{weakly in} \quad X, \quad i=1,2, \dots, N. \label{P_weak}
 \end{align}
  Next we show that  $P^*$ is the solution corresponding to  $U^*$.
According to Definition \ref{def1}, 
  \begin{align}
( \dot{P^k}(t),Q)=(\A P^k(t), Q)+( \sum^{N}_{i=0} f^k(t-t_{i}), Q), \quad \forall Q\in X^*,  \label{weak_form1}
 \end{align}
 and $P^k(0)=P_0$.  Note that 
 \[(\A P, Q\psi)=(P, \A^*Q\psi),\]
 where  $A^*$ is the adjoint operator of $A$ defined by
   \begin{eqnarray}
\A^* Q=
    \left[\begin{array}{c}
-\lambda_{1}  (q_{0}-q_{1}(0))\\
 \frac{dq_{1}(x)}{dx}+\mu_{1}(x)(q_{0}-q_{1}(x))
\end{array}\right], 
\end{eqnarray} 
with its domain 
\begin{eqnarray*}
D(\A^*)=\big\{Q\in X^*\big| \ \frac{ d q_1} {d x} \in L^{\infty}[0, l],  \mu_1q_{1}\in L^{\infty}[0, l],  \ \text{and}\ q_1(l)=0 \big\}.
 \end{eqnarray*}
Let $\psi=[\psi_0, \psi_1]^T$ be a vector of  continuously differential functions on $[0, T]$ with $\psi_j(T)=0, j=0,1$. For each $Q\in X^*$, we multiply \eqref{weak_form1} by $\psi$ and integrate the  left hand side of  \eqref{weak_form1} by parts with respect to $t$.
 This process yields 
   \begin{align}
-\int^T_{0}( P^k(t),Q\dot{\psi})\, dt=&\int^T_{0}(P^k(t), \A ^*Q\psi)\, dt \nonumber\\
& +\sum^N_{i=1}\left(\int^l_{0}b^k_{i}(x)p^k_1(x, t_i)\, dx\right)q_0(t_i) \psi_{0}(t_i)\nonumber\\
&-\sum^{N}_{i=0} (b^k_i(x)p^k_1(x, t_i), q_1(x, t_i)\psi(t_i))+(P^k(0),Q\psi(0)), \label{weak_form2}
 \end{align}
 where by virtue of \eqref{U_uniform}--\eqref{P_weak}, we have
\begin{align}
& \left|\int^l_{0}b^k_{i}(x)p^k_1(x, t_i)\, dx-\int^l_{0}b^*_{i}(x)p^*_1(x, t_i)\, dx\right|\nonumber\\
&= \left|\int^l_{0}(b^k_{i}(x)-b^*(x))p^k_1(x, t_i)\, dx- \int^l_{0}b^*_{i}(x)(p^*_1(x, t_i)-p^k_1(x, t_i))\, dx\right|\nonumber\\
&\leq    \|b^k_{i}-b^*\|_{L^{\infty}[0, l]}\|p^k_1(\cdot, t_i)\|_{L^1[0, l]}+ \left|\int^l_{0}b^*_{i}(x)(p^*_1(x, t_i)-p^k_1(x, t_i))\, dx\right|\to 0.
  \label{EST_bq}
 \end{align}
 Moreover, 
 \begin{align}
&\left|\left(b^k_i(x)p^k_1(x, t_i), q_1(x, t_i)\right)- \left(b^*_i(x)p^*_1(x, t_i), q_1(x, t_i)\right) \right|\nonumber\\
 &=\left | \left((b^k_i(x)-b^*_i(x))p^k_1(x, t_i), q_1(x, t_i)\right)+
 \left(b^*_i(x)(p^k_1(x, t_i)-p^*_1(x, t_i)), q_1(x, t_i)\right)\right| \nonumber\\
& \leq \|b^k_i(x)-b^*_i(x)\|_{L^{\infty}[0, l]}\|p^k_1(x, t_i)\|_{L^{1}[0, l]}\|q_1(\cdot, t_i)\|_{L^{\infty}[0, l]} \nonumber\\
&\quad +\left|\left(p^k_1(x, t_i)-p^*_1(x, t_i), b^*_i(x)q_1(x, t_i)\right)\right|\to 0, \label{EST_bpq}
 \end{align}
 where we used the fact  that $b^*_i(x)q_1(x, t_i)\in L^{\infty}[0, l]$  for the second term of \eqref{EST_bpq} converging to zero.
 
Now  pass to the limit in each term of \eqref{weak_form1} by using \eqref{weak_form2}--\eqref{EST_bpq}. As a result, we have
    \begin{align}
-\int^T_{0}( P^*(t),Q\dot{\psi})\, dt=&\int^T_{0}(P^*(t), \A^*Q\psi)\, dt\nonumber\\
&+ \int^T_{0}( \sum^{N}_{i=0} f^*(t-t_{i}), Q\psi)\, dt+(P_0,Q\psi(0)). \label{weak_form3}
 \end{align}
 where $ f^*(t-t_{i})=[\delta(t-t_i)\int^l_{0}b^*_i(x)p^*_{1}(x,t_i)\,dx, -\delta(t-t_i)b_i(x)p_{1}(x,t_i)]^T $.
  It remains to  be shown that $P^*(0)=P_0$.
  Consider
  \begin{align}
( \dot{P}^*(t),Q)&=(\A P^*(t), Q)+( \sum^{N}_{i=0} f^*(t-t_{i}), Q), \quad \forall Q=[q_0, q_1]^T\in X^* \label{weak_Pstar}\\
P^*(0)&=P^*_0. \nonumber
 \end{align}
 We repeat  the same process as above. Multiplying \eqref{weak_Pstar}   by a continuously differentiable function  $\psi$  with $\psi(T)=0$ and integrating  by parts yield 
    \begin{align}
-\int^T_{0}( P^*(t),Q\dot{\psi})\, dt=&\int^T_{0}(P^*(t), \A^*Q\psi)\, dt+ \int^T_{0}( \sum^{N}_{i=0} f^*(t-t_{i}), Q\psi)\, dt\nonumber\\
&+(P^*_0,Q\psi(0)). \label{weak_form4}
 \end{align}
 Comparing \eqref{weak_form4} with \eqref{weak_form3} gives
 \begin{align}
 (P^*_0,Q\psi(0))=(P_0,Q\psi(0)), \quad \forall Q\in X^* \label{EST_ini}
 \end{align}
Choose $\psi$ with $\psi(0)=1$. Then \eqref{EST_ini} becomes 
\[(P^*_0-P_{0}, Q)=0, \quad \quad \forall Q\in X^*,\]
and thus $P^*_0=P_{0}$.
Lastly, by the  lower semicontinuity of $J$  for all $U\in U_{ad}$, we have
\begin{align*}
J(U^*)\leq \lim_{k\to\infty}\inf J(U^{k}).
\end{align*}
 This completes the proof. 
 \section {Optimality Conditions}
In this section, the first-order  necessary  optimality conditions  for problem \eqref{cost} will derived  by using a variational inequality \cite{lion1971}.
If $U$ is an optimal solution of problem \eqref{cost}, then
\begin{align} 
J'(U)\cdot (V-U)\geq 0, \quad  \forall V\in U_{ad},\label{var_ineq}
\end{align}
where $J'(U)\cdot h$ stands for  the G$\hat{a}$teaux derivative of $J$ at $U$ in the direction  $h\in U_{ad}$. 

Define  operator $D\colon C([0, T], L^{1}[0, l])\to C[0, T] $ by
 \[D p_{1}(x,t)=\int^l_{0}p_1(x,t)\, dx=\hat{p}_1(t).\]
According to the definition of $J$ in \eqref{cost}, we have
\begin{align}
J'_i(U)\cdot h_i  =&\int^{T}_{0}(D^*D p_{1}(x,t), z_{1i})\,d t
+\beta_i \int^l_{0}b_{i}h_i\, dx+ (D^*Dp_{1}(x, T), z_{1i}(x, T) ), 
 \label{var_J}
\end{align} 
for $i=1,2\dots, N$, where $z_{0i}=p'_{0}(u_{i})\cdot h_{i}$, $z_{1i}=p'_{1}(u_{i})\cdot h_i$, and $h=[h_1, h_2, \dots, h_N]^T\in U_{ad}$.  Note that   the G$\hat{a}$teaux derivatives $z_{0i}$ and $z_{1i}$ satisfy the following equations 
 \begin{align}
&\frac{dz_{0i}(t)}{d t}=-\lambda_1  z_{0i}(t)+ \int_0^{l} \mu_1(x)z_{1i}(x,t) \,d x  \nonumber\\
&\qquad+  \delta(t-t_{i}) \int_0^{l}(h_{i}(x)p_1(x, t_i)+b_{i}(x)z_{1i}(x, t_i)\,d x,  \label{var_1}\\
&\frac{\partial{z_{1i}(x,t)}}{\partial{t}}+\frac{\partial{z_{1i}(x,t)}}{\partial{x}} = -\mu_1(x)z_{1i}(x,t) \nonumber\\
&\qquad -\delta(t-t_{i}) (h_{i}(x)p_1(x, t_i)+b_{i}(x)z_{1i}(x, t_i)),
 \label{var_2}
\end{align}
with boundary conditions
\begin{equation}
z_{1i}(0,t)=\lambda_1 z_{0i}(t), \label{var_BC}
\end{equation}
and initial conditions
\begin{equation}
z_{0i}(0)=0, \quad z_{1i}(x,0)=0, \quad i=1,2\dots, N. \label{var_INI}
\end{equation}
The well-posedness of  \eqref{var_1}--\eqref{var_INI} can be established by using the similar  approach as  in Theorem \ref{T01}.

Finally, the  first-order necessary conditions of optimality are given by the following theorem.
\begin{theorem}
Assume that $U^*=[b^*_1, b^*_2, \dots, b^*_N]^T\in U_{ad}$ is  an  optimal solution of   \eqref{cost} subject to \eqref{IVP_A}--\eqref{IVP_IC2}.  Let  $[p_0, p_1]^T\in PWC_l([0, T]; X)$ be the corresponding solution and $[q_{0i},q_{1i}]^T\in C[0, T]\times PWC_r([0, T]; L^{\infty}[0, l]), i=1,2, \dots, N$,  be the solutions to  the following
adjoint systems 
  \begin{align}
&-\frac{dq_{0i}}{dt}=-\lambda_{1}  (q_{0i}-q_{1i}(0)), \quad q_{0i}(T)=0, \label{adj_1} \\
& -\frac{\partial{q}_{1i}}{\partial{t}} -\frac{\partial{q}_{1i}}{\partial{x}}=\mu_{1}(x)(q_{0i}(t)-q_{1i}(x, t))+D^*Dp_{1}(x,t)\nonumber\\
&\qquad + \delta(t-t_i)b_{i}(x) q_{1i}(x, t_i) -  \delta(t-t_i)\int_0^{l}b_{i}(x)\,d x \,q_{0i}(t_i),\label{adj_2}\\
 &  q_{1i}(l, t)=0, \quad q_{1i}(x, T)=D^*Dp_{1}(x, T)=\int^{l}_{0}p_{1}(x, T)\, dx, \quad i=1,2, \dots, N.
\label{adj_3}
\end{align}
Then   $b^*_i$  satisfies   
 \begin{align}
 b^*_i=\max\{0, \min\{ \beta^{-1}_i p_1(x, t_i)(q_1(x, t_i)- q_0(t_i) ), 1\}\}, \quad i=1,2,\dots, N,
 \end{align}
  where $PWC_r([0, T]; L^{\infty}[0, l])$ denotes the space of piecewise continuous functions  on $[0, T]$ with values in $L^{\infty}[0, l]$, that are right continuous and possess lefthand limits.

\end{theorem}
 {\noindent {\textbf{Proof.} } 
We first compute  the first term in  \eqref{var_J}. This yields
 \begin{align}
&  \int^{T}_{0}(D^*D p_{1}(x,t), z_{1i})\,d t=\int^{T}_{0}\bigg( -\frac{\partial{q}_{1i}}{\partial{t}} -\frac{\partial{q}_{1i}}{\partial{x}}-\mu_{1}(x)(q_{0i}(t)-q_{1i}(x, t))\nonumber\\
&\qquad+  \delta(t-t_i)  b_{i}(x) q_{1i}(x, t_i) - \delta(t-t_i)  \int_0^{l}b_{i}(x)\,d x q_{0i}(t_i), z_{1i}\bigg)\,d t
 \nonumber\\
&=\int^l_{0}(-q_{1i}(x, T)z_{1i}(x, T)+q_{1i}(x, 0)z_{1i}(x, 0))\, dx+\int^{T}_{0} (\frac{\partial{z}_{1i}}{\partial{t}}, q_{1i})\, dt\nonumber\\
&\qquad+\int^T_{0}(-q_{1i}(l, t)z_{1i}(l,t)+q_{1i}(0, t)z_{1i}(0, t))\, dt+\int^{T}_{0} (\frac{\partial{z}_{1i}}{\partial{x}},q_{1i})\, dt\nonumber\\
&\qquad+ \int^{T}_{0}\bigg( -\mu_{1}(x)(q_{0i}(t)-q_{1i}(x, t))+\delta(t-t_i)  b_{i}(x) q_{1i}(x, t_i) \nonumber\\
&\qquad -    \delta(t-t_i)   \int_0^{l}b_{i}(x)\,d x q_{0i}(t_i), z_{1i}\bigg)\,d t\nonumber\\
&=\int^l_{0}(-q_{1i}(x, T)z_{1i}(x, T)\, dx+\int^{T}_{0}( \frac{\partial{z}_{1i}}{\partial{t}}, q_{1i})\, dt
+\int^T_{0}(q_{1i}(0, t)\lambda_1z_{0i}( t))\, dt\nonumber\\
&\qquad +\int^{T}_{0} (\frac{\partial{z}_{1i}}{\partial{x}}, q_{1i})\, dt
+ \int^{T}_{0}\bigg( -\mu_{1}(x)(q_{0i}(t)-q_{1i}(x, t))+  \delta(t-t_i)  b_{i}(x) q_{1i}(x, t_i)\nonumber\\
&\qquad-  \delta(t-t_i)   \int_0^{l}b_{i}(x)\,d x q_{0i}(t_i), z_{1i}\bigg)\,d t. \label{J_first}
\end{align}
Note that   by  \eqref{var_2} and \eqref{adj_1} we have
\begin{align}
&\int^{T}_{0}( \frac{\partial{z}_{1i}}{\partial{t}}, q_{1})\, dt+\int^{T}_{0} (\frac{\partial{z}_{1i}}{\partial{x}}, q_{1})\, dt \nonumber\\
&\qquad=\int^T_{0} \bigg(-\mu_1(x)z_{1i}(x,t) 
- \delta(t-t_{i}) (h_{i}(x)p_1(x, t_i)+b_{i}(x)z_{1i}(x, t_i)), q_1\bigg)\, dt \label{EST_z}
\end{align}
and 
\begin{align}
\int^T_{0}(q_{1}(0, t)\lambda_1z_0( t))\, dt=\int^T_{0}(-\frac{dq_{0}}{dt}+\lambda_{1}  q_{0})z_0\, dt. \label{EST_q0}
\end{align}
Thus combining \eqref{J_first} with \eqref{EST_z}--\eqref{EST_q0} yields
\begin{align}
&  \int^{T}_{0}(D^*D p_{1}(x,t), z_{1i})\,d t=-\int^l_{0}q_{1i}(x, T)z_{1i}(x, T)\, dx +\int^T_{0}(-\frac{dq_0}{dt}z_{0i}+\lambda_1q_{0i}z_{0i})\,  dt
\nonumber\\
&\quad+\int^T_{0} \bigg(-\mu_1(x)z_{1i}(x,t) 
-\delta(t-t_{i}) (h_{i}(x)p_1(x, t_i)+b_{i}(x)z_{1i}(x, t_i)), q_{1i}\bigg)\, dt  \nonumber\\
&\quad+ \int^{T}_{0}\bigg( -\mu_{1}(x)(q_{0i}(t)-q_{1i}(x, t))+ \delta(t-t_i)  b_{i}(x) q_{1i}(x, t_i) \nonumber\\
&\quad-   \delta(t-t_i)    \int_0^{l}b_{i}(x)\,d x q_{0i}(t_i), z_{1i}\bigg)\,d t\nonumber\\
=&-\int^l_{0}q_{1i}(x, T)z_{1i}(x, T)\, dx +(-q_{0i}(T)z_{0i}(T)+q_{0i}(0)z_{0i}(0))\nonumber\\
&+\int^T_{0}(\frac{dz_{0i}}{dt}q_{0i}+\lambda_1q_{0i}z_{0i})\,dt+\int^T_{0} \bigg(-\delta(t-t_{i}) h_{i}(x)p_1(x, t_i), q_{1i}\bigg)\, dt \nonumber\\
&+ \int^{T}_{0}\bigg( -\mu_{1}(x)q_{0i}(t)- \delta(t-t_i)    \int_0^{l}b_{i}(x)\,d x q_{0i}(t_i), z_{1i}\bigg)\,d t\nonumber\\
=&-\int^l_{0}q_{1i}(x, T)z_{1i}(x, T)\, dx \nonumber\\
&+ \int^T_{0}\bigg(\int_0^{l} \mu_1(x)z_{1i}(x,t) \,d x 
+ \delta(t-t_{i}) \int_0^{l}(h_{i}(x)p_1(x, t_i)+b_{i}(x)z_{1i}(x, t_i)\,d x\bigg)q_{0i}\, dt\nonumber\\
&+\int^T_{0} \bigg(- \delta(t-t_{i}) h_{i}(x)p_1(x, t_i), q_{1i}\bigg)\, dt  \nonumber\\
&+ \int^{T}_{0}\bigg( -\mu_{1}(x)q_{0i}(t)-   \delta(t-t_i)    \int_0^{l}b_{i}(x)\,d x q_{0i}(t_i), z_{1i}\bigg)\,d t\nonumber\\
=&-\int^l_{0}q_{1i}(x, T)z_{1i}(x, T)\, dx + \int^T_{0}\bigg( \delta(t-t_{i}) \int^l_0h_{i}(x)p_1(x, t_i)\,d x\bigg)q_{0i}\, dt\nonumber\\
&+\int^T_{0} \bigg(- \delta(t-t_{i}) h_{i}(x)p_1(x, t_i), q_{1i}\bigg)\, dt \nonumber\\
=&-\int^l_{0}q_{1i}(x, T)z_{1i}(x, T)\, dx +  \int^l_0h_{i}(x)p_1(x, t_i)\,d x\, q_{0i}(t_i) \nonumber\\
&-\int^l_{0}h_{i}(x)p_1(x, t_i)q_{1i}(x, t_i)\, dx. \label{EST_J}
 \end{align}
Replacing   $\int^{T}_{0}(D^*D p_{1}(x,t), z_{1i})\,d t$ in \eqref{var_J} by   \eqref{EST_J} and making use of the condition \eqref{adj_3}, we have  the  G$\hat{a}$teaux  derivative $J'_{i}(U)\cdot h_i$ become
 \begin{align}
J'_i(U)\cdot h_i  =&-\int^l_{0}q_{1i}(x, T)z_{1i}(x, T)\, dx +  \int^l_0h_{i}(x)p_1(x, t_i)\,d x\, q_{0i}(t_i)\nonumber\\
&-\int^l_{0}h_{i}(x)p_1(x, t_i)q_{1i}(x, t_i)\, dx
+\beta_i \int^l_{0}b_{i}h_i\, dx\nonumber\\
&+ (D^*Dp_{1}(x, T), z_{1i}(x, T) )\nonumber\\
=&    \int^l_0h_{i}(x)p_{1i}(x, t_i)(q_{0i}(t_i)
-q_{1i}(x, t_i))\, dx+\beta_i \int^l_{0}b_{i}h_i\, dx \geq 0, \label{opt_cond1}
\end{align}
for $i=1,2, \dots, N,$ and any $h=[h_1, h_2, \dots, h_N]^T\in U_{ad}$. Finally, combining \eqref{opt_cond1} with the constraint  that $0 \leq b_i\leq 1$, we get 
 \begin{align*}
 b_i=\max\{0, \min\{1, \beta^{-1}_i p_1(x, t_i)(q_{1i}(x, t_i)- q_{0i}(t_i) )\}\},  \quad i=1,2, \dots, N.
 \end{align*}
 This completes the proof.

\section{Conclusion}

An  impulse control design  is discussed for a simple reparable system in a nonreflexive Banach space, which represents a  preventive maintenance policy. 
First-order conditions of optimality are derived for solving the optimal solution.  A finite difference scheme that preserves the nonnegativity and  conservativeness   will be developed to discretize  the controlled system and  a gradient decent based algorithm will be constructed  to implement the control design in our future work.

\section*{Acknowledgments} 
H. Xu was supported by  NSAF grant No. U 1430125.








\end{document}